\documentclass[journal]{IEEEtran}

\ifCLASSINFOpdf
\else
\fi
\usepackage{times,amsmath,epsfig}
\usepackage{amssymb, dsfont}
\usepackage{booktabs}
\usepackage{breqn}
\usepackage[compress]{cite}
\usepackage{tabularx}
\usepackage{amsfonts}
\usepackage{algorithm}
\usepackage{algorithmic}
\usepackage{caption}
\usepackage[labelformat=simple,position=b]{subcaption}
\usepackage{tikz}
 \usepackage[usenames,dvipsnames]{pstricks}
\usetikzlibrary{automata,positioning}
\usetikzlibrary{arrows}
\usetikzlibrary{calc}
\usepackage{comment}
\usepackage{flushend}
\usepackage[intoc]{nomencl}
\newcommand{\BigO}[1]{\ensuremath{\operatorname{O}\bigl(#1\bigr)}}
\makenomenclature

\makeindex 

\definecolor {processblue}{cmyk}{0,0,0,0.17}
\definecolor{light-gray}{gray}{0.95}

\begin{document}

\title{Capacity Planning Frameworks for Electric Vehicle Charging Stations with Multi-Class Customers}
\author{Islam~Safak~Bayram,~\IEEEmembership{Member,~IEEE,}~Ali~Tajer,~\IEEEmembership{Member,~IEEE,}~Mohamed~Abdallah,~\IEEEmembership{Senior~Member,~IEEE}~and
~Khalid~Qaraqe,~\IEEEmembership{Senior~Member,~IEEE}

\thanks{Islam Safak Bayram (corresponding author) is with Qatar Environment and Energy Research Institute, Qatar Foundation, Doha Qatar. Email: ibayram@qf.org.qa}
 \thanks{Mohamed Abdallah and Khalid Qaraqe are with the Department
of Electrical and Computer Engineering, Texas A\&M University at Qatar. Emails:\{mohamed.abdallah, khalid.qaraqe\}@qatar.tamu.edu.}
\thanks{Ali Tajer is with the Department of Electrical, Computer, and Systems Engineering, Rensselaer Polytechnic Institute, USA, Email: tajer@ecse.rpi.edu.}}

\markboth{{IEEE Transactions on Smart Grid \emph{accepted for publication}}}%
{Shell \MakeLowercase{\textit{et al.}}: Bare Demo of IEEEtran.cls for Journals}

\maketitle
\IEEEpeerreviewmaketitle

\begin{abstract}
In order to foster electric vehicle (EV) adoption, there is a strong need for designing and developing charging stations that can accommodate different customer classes, distinguished by their charging preferences, needs, and technologies. By growing such charging station networks, the power grid becomes more congested and, therefore, controlling of charging demands should be carefully aligned with the available resources. This paper focuses on an EV charging network equipped with different charging technologies and proposes two frameworks. In the first framework, appropriate for large networks, the EV population is expected to constitute a sizable portion of the light duty fleets. This which necessitates controlling the EV charging operations to prevent potential grid failures and distribute the resources efficiently. This framework leverages pricing dynamics in order to control the EV customer request rates and to provide a charging service with the best level of quality of service. The second framework, on the other hand, is more appropriate for smaller networks, in which the objective is to compute the minimum amount of resources required to provide certain levels of quality of service to each class. The results show that the proposed frameworks ensure grid reliability and lead to significant savings in capacity planning.
\end{abstract}

\printnomenclature
\section{Introduction}

\subsection{Motivation}
Electric Vehicles are becoming a viable transportation option as they offer solutions to an array of current societal problems ranging from high oil prices to environmental concerns. Corollary,  EV population is expected to reach a sizable market portion in the next decade (e.g., $10$\% of the U.S. National fleet and similar targets in Europe)~\cite{jsac}. However, achieving such penetration rates requires wide deployment of charging facilities that can serve different types of charging requests. On the other hand, if not controlled, EV charging can easily lead to transformer and line overloading. Moreover it can deteriorate the power quality and even endangers the security of supply. The impacts of EV charging on the grid has been very well documented in~\cite{5356176, jsac,generationPortfolio} and \cite{oakRidge}. For instance \cite{oakRidge} argues that if $5\%$ of the EVs charge simultaneously using fast charging technology,  $5.5$ GW of extra power would be needed in Virginia and Carolinas region by $2018$.  Similarly, NERC regions would require an expansion of $5.5$\% in their power generation capacity in a typical EV penetration rate of $25$\%. Furthermore distribution grid could easily become a bottleneck. For example, authors of \cite{jsac} state that even adding two Level-2 chargers in a typical neighborhood in the US could easily cause transformer overloading. Also uncontrolled demand can decrease the efficiency of the power grid operations and increase the power generation cost if it occurs during peak hours.
\subsection{Related Works}\label{relatedWork2}

There has been an increasing body of literature on developing intelligent charging station architectures\cite{power1,power2,power4,sgc12, queue2, queue1,queue3}, controlling and scheduling EV demand \cite{TSG14,jsac,CallawayX,price1,price2,price3,globecom} and rather sparse literature on capacity planning on charging infrastructures serving multiple classes of customers \cite{sgc13,energyCon}. This section provides a brief overview of related literature.

The studies on charging station design can be divided into two categories. The first category includes works from power engineering perspective \cite{power1,power2,power4} where the studies aim to minimize the charging duration by improving the efficiency of power electronics and aiding the system with energy storage system. On the other hand, approaches in the second category are mostly focused on the system level where they abstract the details of the underlying power system components and evaluate the system performance in terms of long term statistical metrics e.g., mean waiting time and percentage of served customer using the arguments from queuing theory~\cite{ sgc12, queue2, queue3}. Study in \cite{sgc12} presents a small scale fast charging station design and blocking probability, that is the probability that the request of an EV is rejected, is used as the main performance metric to solve the optimal resource provisioning problem. In\cite{queue1} authors use M/M/s queues to model the EV demand at fast charging stations near highway exits. In~\cite{queue2} residential charging infrastructures are modeled using M/M/$\infty$ queue, where they broadcast the arrival rates to customers such that probability of blocking (exceeding circuit capacity) is zero. Furthermore, the work in~\cite{queue3} uses BCMP queueing network model to estimate the EV charging demand. In this paper we perform a more holistic approach by taking into account that the charging station can serve multiple classes of customer that are differentiated by the charger technology used. Also, we employ loss-of-load-probability (LoLP) as our main performance metric, which measures the probability that grid resources cannot accommodate EV demand \cite{von2006electric}.

The pricing-based control methods have been employed successfully in communication networks where the goal is to match scarce bandwidth resources to insatiable user demand~\cite{priceSurvey}. In a similar fashion, pricing mechanisms applications gained popularity in the literature on controlling EV demand~\cite{jsac,globecom, price1, price2, price3} and \cite{fan2011}. Authors in \cite{jsac} and \cite{globecom} propose a control mechanism in a network of charging stations to route customers to idle stations. The work in~\cite{price1} proposes a pricing scheme for EV chargings that leads to socially optimal solution, whereas \cite{fan2011} uses proportional fairness pricing from communication literature to control the charging rates of EVs. In our work we also consider the technological constraints and current state of affairs that charging rate is constant until the end of the service. The work presented in \cite{price3} is most relevant to our work, in which researchers model general power consumers behavior with utility functions and propose a pricing algorithm to control the consumption of smart grid users. Our work mainly focuses on multi-class EV chargings and our framework ensures grid reliability at all times. The literature on resource provisioning is rather sparse. The study in \cite{sgc13} presents a capacity planning framework in a large scale charging stations with single class customers. This framework is based on computing a deterministic quantity ``effective power" using On-Off Markov models. In this study we assume that charging infrastructure can serve multiple types of customers.

The proposed methodology is rooted in multi dimensional loss systems in teletraffic engineering, where the goal is to provide statistical quality of service(QoS) guarantees to customers with different demand profiles. The control mechanisms with  congestion pricing in multi-rate Erlang-B systems and related resource provisioning problems are addressed in~\cite{hampshire, kaufman, nilsson, cong1,cong2}. The works presented in \cite{cong1} and \cite{kaufman,nilsson,fan2011} focus on efficient computation of blocking probabilities and the their derivatives, whereas studies in \cite{hampshire} and \cite{cong2} provide efficient algorithms to solve bandwidth provisioning problem in congestible networks.
\begin{figure}
                \centering
                \includegraphics[width=0.7\columnwidth]{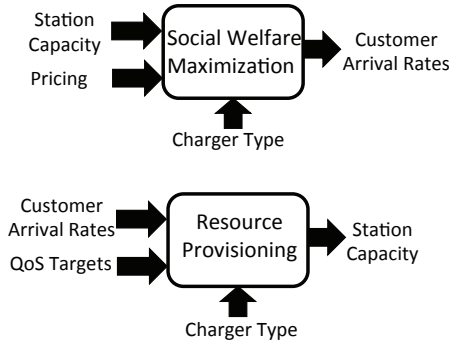}
 \caption{Two design problems. }\label{comp1}
        \end{figure}

\section{Proposed Frameworks}\label{ChargingStationModel}

Electric vehicle charging demands are primarily dominated by customers' needs and preferences. Based on this premise, designing different aspects of the network of charging stations is governed by the  economic dynamics between network operators, on one hand, and the customers on the other hand. In such dynamics the former seek operational reliability and profit margins and the latter seek economic incentives. Such dynamics and their pertinent design challenges vary by the size of network and customers population. In this paper we consider two canonical design problems for EV charging infrastructures fed by a single substation suited for large metropolitan areas and small cities as  summarized in Fig.~\ref{comp1}.

Specifically, the first framework proposes a pricing-based EV control for charging stations located in large metropolitan cities in which the EV population is assumed to be large. In this framework, the objective is to incentivize and control the customers to submit their charging demands at certain optimal rates which maximize a social welfare utility. This pricing-based approach framework, at its core, enables station operators to alleviate the congestion and the degradations in the service quality during increased demand periods through incentivizing the customers to defer their charging needs to less congested periods. This control mechanism ensures charging services with desirable level of QoS guarantees.

The second framework focuses on small cities, in which customer demand can be obtained from profiling studies, and proposes a capacity provisioning mechanism for EV charging stations. In contrary to large metropolitan areas, small cities with fewer charging stations and customers enjoy less flexibility for incentivizing the customers. Hence, the primary goal in this framework is to compute the minimum amount of required resources to provide charging service to meet a target level of QoS.

\subsection{Network Model}\label{sysModel}
\nomenclature{$N$}{Number of customers}%
\nomenclature{$i$}{Customer index }%
\nomenclature{$J$}{Number of distinct customer classes}%
\nomenclature{$j$}{Customer class index}%
\nomenclature{$b_j$}{Amount of resources requested by customer class $j$}%
\nomenclature{$k$}{Time index}%
\nomenclature{$C_k$}{Amount of power drawn from the grid during period $k$}%
\nomenclature{$\beta_j(\cdot)$}{Loss-of-Load-Probability for class $j$}%
\nomenclature{$p_j(k)$}{Price of service for customer type $j$ at time $k$}%
\nomenclature{$\lambda_j^k$}{Arrival rate for customer class $j$ at time $k$}%
\nomenclature{$\lambda_j(\cdot)$}{Aggregated arrival rate for customer class $j$ at time $k$ given in \eqref{aggrArrival}}%
\nomenclature{${\boldsymbol \lambda}^n$}{Arrival rate vector for all customer classes at period $n$}%
\nomenclature{${\boldsymbol \lambda}$}{Aggregated arrival rate vector for all customer classes at period $n$}%
\nomenclature{$U^n(\cdot)$}{Utility function for customer $n$}%
\nomenclature{$\delta_j$}{Loss-of-Load Probability target for customer class $j$}%
We consider a network of charging stations serving $N$ customers. Since the customers do not necessarily seek identical charging services, we consider $J$ distinct customer classes, which are distinguished by their preferences; size of the battery packs, amounts of requested demands, and the available charger technologies. The charging rate of customers of type $j \in \left\{ {1,..,J} \right\}$ is denoted by $b_j$. Since the amount of resources are finite, upon the arrival of a specific customer type $j$, if the available resources is less than $b_j$ the customer will not be served resulting in an outage. Hence, the probability of being blocked (or loss-of-load-probability) constitutes a natural system performance metric. Due to the temporal fluctuations network conditions (e.g., congestion etc.), the amount of available resources varies over time. Hence, we consider a dynamic system indexed by the time index $k\in\mathds{N}$ in which $C_k$ denotes the aggregate units of grid power available to the entire EV fleet at time $k$. Note that $C_k$ is an exogenous parameter given to station. Since the charging stations reside in a small well-confined regions, the set of chargers can be abstracted collectively as a super-station with multiple classes of customers. We further assume that system capacity is large when compared to demands of each EV type.

Increased outage events leads to reduced utility for the network operator and service disruption for the customers. Hence, optimizing the operation of the network strongly hinges on finding the decision rules that lead to the optimal level of Loss-of-Load-Probabilities (LoLP) that is defined as the probability that grid resources fall short of aggregated customer demand \cite{von2006electric}. The leverage that the network operators have for adjusting the LoLP rates at the desired optimal rates is service pricing, by using which they can influence the behavior of the customers abstracted by their arrival rates. To formalize this, we define $p_j(k)$ as the price of service of type $j$ and define the price vector $\boldsymbol{p}(k)\triangleq[p_1(k),\dots,p_J(k)]$, accordingly.
 \tikzstyle{block} = [rectangle, draw,top color =light-gray , bottom color = processblue!20,
    text width=5em, text centered, rounded corners, minimum height=3em]
\tikzstyle{line} = [draw, -triangle 45]

\begin{figure}
\begin{center}
\begin{tikzpicture}[node distance = 2.5cm, auto,scale=1, every node/.style={transform shape}]
\node [block] (init) {Aggregated Load};
\node [block, left of=init] (Grid) {Grid Conditions};
\node [block, below of=Grid] (Demand) {EV Demand};
\node [block, right of=Demand] (EVAC) {Pricing};
\path [line] (Grid) -- (init);
\path [line,dashed] (Demand) -- (init);
\path [line,dashed] (init) -- (EVAC);
\path [line,dashed] (EVAC) -- (Demand);
\end{tikzpicture}
 \caption{Pricing-based control mechanism for social welfare maximization }\label{illustrate}
 \end{center}
 \end{figure}
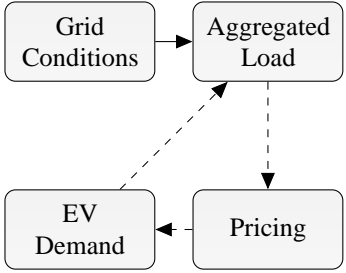
Based on the prices ${\boldsymbol p}(k)$ and its needs, each customer decides whether to generate a service request of a certain type. We denote the rate of  type $j \in \left\{ {1,..,J} \right\}$ requests generated by customer $n\in \left\{ {1,..,N} \right\}$ during period $k$ by ${\lambda _j^{n}(k;{\boldsymbol p}(k))}$. Hence, the aggregate arrival rate of service requests of type $j$ is
\begin{equation}\label{aggrArrival}
\lambda_j(k;{\boldsymbol p}(k))\triangleq\sum_{n=1}^N\lambda_j^n(k;{\boldsymbol p}(k)) \ .
\end{equation}
From the network operators point of view, these requests are being constantly generated across the networks and as discussed in \cite{queue1,queue2,jsac}, a common assumption is that such aggregated requests (and not necessarily individual requests) arrive according to a Poisson process. Accordingly, we define the arrival rate vectors
\begin{align}
{\boldsymbol \lambda}^n(k;{\boldsymbol p}(k))\; & \triangleq \; [\lambda^n_1(k;{\boldsymbol p}(k)),\dots,\lambda^n_J(k;{\boldsymbol p}(k))]\ ,\\
\mbox{and}\qquad {\boldsymbol \lambda}(k;{\boldsymbol p}(k))\; & \triangleq[\lambda_1(k;{\boldsymbol p}(k)),\dots,\lambda_J(k;{\boldsymbol p}(k))]\ .
\end{align}
By noting that the LoLPs for each class are functions of the arrival rates, finally we define $\beta_j(k;{\boldsymbol \lambda}(k;{\boldsymbol p}(k))$ as the probability that customers of type $j$ are blocked and define ${\boldsymbol \beta}(k;{\boldsymbol \lambda}(k;{\boldsymbol p}(k))$ as the corresponding LoLP vector.

Since different customer types do not necessarily charge their batteries at the same rate, we denote the {\em average} charging duration of customers of type $j$ by $1/{\mu _j}$. It is noteworthy that from station operator standpoint, the information required to characterize customers' profiles are the average service duration and the charging rates. These parameters enable modeling customers' behaviors with fine granularity.

The main goal is to keep the aggregated demand below a constant power $C_{k}$ with minimal or controlled loss-of-load (outage) events. Drawing constant power is highly desirable from the viewpoint of the network operators and its benefits include (1) grid components are isolated from stochastic variations and hence grid reliability is ensured~\cite{sgc12,jsac}; (2) station operator can sign long-term contracts and benefit from the lower prices; (3) constant demand will reduce the peak-to-average demand ratio of the whole power system and accordingly the average spot prices would reduce; and (4) it leads to more efficient market equilibrium \cite{jsac}. It is noteworthy that since the charging stations reside in a small well-confined region, the whole set of chargers act as one big station with multiple classes of customers. Due to this fact, power system losses are assumed to be negligible and all customers share the same resource pool.

For the simplicity of the notations, in the rest of the paper the explicit dependence on $k$ is omitted.
\tikzstyle{block} = [rectangle, draw,top color =light-gray , bottom color = processblue!20,
    text width=5em, text centered, rounded corners, minimum height=3em]
\tikzstyle{line} = [draw, -triangle 45]

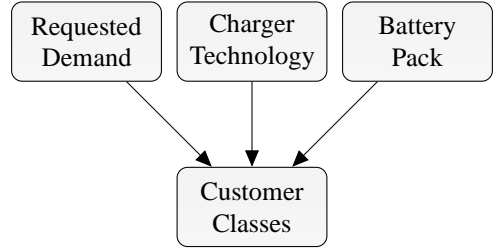
\begin{figure}[t]
\centering
\begin{tikzpicture}[node distance = 2.2cm]
\node [block] (init) {Charger Technology};
\node [block, left of=init] (Grid) {Requested Demand};
\node [block, right of=init] (Policy) {Battery Pack};
\node [block, below of=init] (Demand) {Customer Classes};
\path [line] (Grid) -- (Demand);
\path [line] (init) -- (Demand);
\path [line] (Policy) -- (Demand);
\end{tikzpicture}
\caption{Customer preferences can be determined by EV type, amount of requested demand, and the charger technology.}\label{depend}
\end{figure}

\subsection{Large Networks: Pricing-based Control Framework}\label{pricing}
As the EV population is expected to constitute a sizable portion of the national light duty fleets, there is a need to control EV charging operations in order to protect the grid components from potential failures. To this end, this framework focuses on charging infrastructures that are located in big metropolitan areas with high EV penetration rates. The primary goal is to leverage pricing in order to control the EV customer arrival rates such that station operator can provide a charging service with a small level of customer LoLP QoS for each class. Here, the primary goal is control the large EV demand in a big city for a given amount of grid power ($C_k$), which is assumed to computed by the utility operator and meets the power network constraints (e.g., congestion, transformer ratings etc.).  Pricing-based control scheme is depicted in Fig. \ref{illustrate}.

Specifically, the aim is to design a distributed control framework for maximizing an aggregate utility (social welfare measure), which requires that the EVs operate at an optimal set of arrival rates.  In order to establish the tools for the distributed design let
\begin{align*}
U^n({\boldsymbol \lambda}^{n};{\boldsymbol \beta}({\boldsymbol \lambda}))
\end{align*}
for a given set of arrival rates ${\boldsymbol \lambda}^n$ and block probabilities ${\boldsymbol \beta}$. Hence, the aggregate utility in the network is
\begin{align}\label{mainProblem}
R\;=\; \sum_{n = 1}^N  U^n({\boldsymbol \lambda}^{n};{\boldsymbol \beta}({\boldsymbol \lambda}))\ .
\end{align}
Therefore, the social welfare problem can be defined as the problem that maximizes the aggregate utility $R$ over all possible choices of the arrival rates $\{{\boldsymbol \lambda}^n\}_{n=1}^N$, i.e.,
\begin{align}\label{mainProblem2}
\max_{\{{\boldsymbol \lambda}^1,\dots,{\boldsymbol \lambda}^N\}}  &\sum_{n = 1}^N  U^n({\boldsymbol \lambda}^{n};{\boldsymbol \beta}({\boldsymbol \lambda}))\ .
\end{align}
We assume the utility function $U^n$, is increasing in the arrival rates $\{\lambda_j\}$, and decreasing in the LoLPs $\beta_j$. Furthermore, it is assumed that $U^n$ is concave in $\lambda_j$ and continuously differentiable in all of its arguments. This problem is treated in Section~\ref{sec:price}.

\subsection{Small Networks: Resource Provisioning Framework}\label{resourceProv}
In the resource provisioning problem, we are interested in computing the minimum amount of grid resources $C$ such that system operator guarantees certain levels of reliability for serving different types of users by enforcing $\beta_j \leq \delta_j$ for each customer type. This structure suits the resource provisioning problems in small cities in which customer arrival rates can be obtained from profiling studies with reliable accuracy. Based on the constrains defined for LoLPs, the problem can be cast as
\begin{equation}\label{argmin}
{C^*} = \left\{ {\begin{array}{ll}
\min &  C \\
\mbox{s.t} & \;\beta_j \le {\delta _j},\;\mbox{for}\; \in\{1,\dots,J\}
\end{array}} \right. \ .
\end{equation}
This problem is treated in Section~\ref{lossComp}.
\begin{figure}[t]
 \centering
\includegraphics[width=\columnwidth]{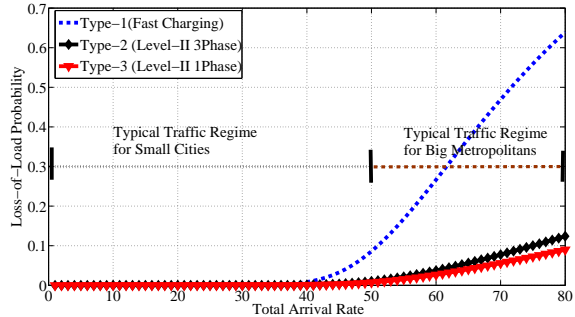}
 \caption{Numerical Evaluation of multi-class EV LoLP. Each EVs arrival rate is assumed to be equal to one third of the total arrival rate.}\label{blockingGraph}
\end{figure}
\subsection{Toy Example}
We provide a toy example before the case studies in order to put some of the details into perspective. Even though the objectives of the proposed frameworks are different, computation of LoLP is carried out for a fixed parameter setting (e.g., $\lambda$, $\mu$ etc.), the methods for computing blocking probability is valid for both cases and details are given in Section \ref{lossComp}. Let us assume that the charging infrastructure draws $C =1000$ units of power from the grid and serves three types of EV customers, i.e.,  $j \in \left\{ {1,2,3} \right\}$. Customer classes are differentiated by the charging technology they use, and by mimicking the current charging standards (fast charging, level-II three and single phase). It is assumed that $\boldsymbol{b}=\left\{ {50,7,5} \right\}$, the charging rates are  $\mu_1=3$, $\mu_2=0.42$, and $\mu_3=0.2$, and the arrival rates for each type are $\lambda_1=\lambda_2=\lambda_3=14$ (on average 14 customers arrive at each time unit). Consequently, the resulting LoLPs are $\beta_{1000}^1$=$0.0152$, $\beta_{1000}^2=0.0015$, and $\beta_{1000}^3=0.0011$. Note that customers of type $1$ have the highest LoLP mainly because grid resources are shared equally among all customer classes and type $1$ customers use fast charging technology which  draws more power compared to other classes.
\begin{figure}[t]
 \centering
\includegraphics[width=\columnwidth]{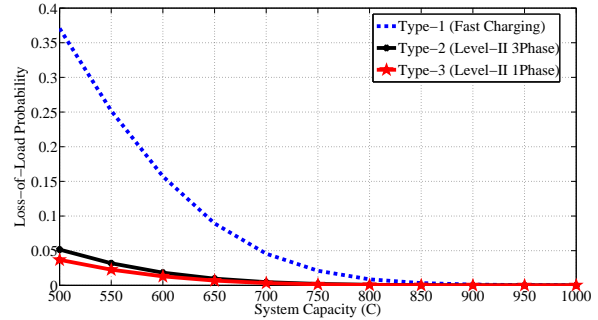}
 \caption{Numerical evaluation of multi-class EV LoLP for different system capacity.}\label{blockingCapacity}
\end{figure}
Next, we evaluate the same system for a wide range of arrival rates. We assume that the sum $\lambda=\sum\nolimits_{j = 1}^3 {{\lambda _j}}$ varies in the range $[1,80]$ and $\lambda_{j}=\lambda / 3$ for $j\in \left\{ {1,2,3} \right\}$. The results are depicted in Fig. \ref{blockingGraph}. It is important to notice that, for the given amount of grid resources (in our case grid allows $C$=$1000$), stations can provide a very good service for light traffic scenarios in which the aggregated arrival rate varies in the range $\lambda\in[1,50]$, which can be a very typical setting for small cities. It is noteworthy that instead of providing an extremely stringent LoLP constraints (near zero outage probability), network operator can back off from $C$=$1000$ to a smaller level $\hat{C}<C$ while guaranteeing a high level of LoLP (e.g. $1$\% LoLP) for each customer classes.

On the other hand, in the heavy traffic regimes (e.g., $\lambda>50$) which can occur in large metropolitan cities, in order to provide good QoS, either charging resources ($C$) should be increased or the arrival rates should be controlled. To shed more light on this, we evaluate the system performance for a fixed set of arrival rates ($\lambda_1$=$\lambda_2$=$\lambda_3=10$) and for a range of serving capacity values $C\in[500,1000]$ . As clearly depicted in Fig. \ref{blockingCapacity}, more customers will be accommodated as the amount of resources increases. Hence, besides upgrading the grid to have increased capacity, which might not be economically viable in the short term, an alternative approach is to provide charging services based on the available resources and with the best possible QoS via controlling customer arrival rates. The pricing-based control framework discussed in section~\ref{pricing} attacks this problem.

\section{Pricing-Based Control}
\label{sec:price}
\label{sec:provision}

\nomenclature{$R$}{Aggregated utility of all customers}%

\subsection{Global Problem}
This treatment provided in this section has a two-fold purpose. Primarily, it provides an optimal LoLP policy in order to maintain an optimal social welfare utility as characterized in \eqref{mainProblem}. Secondly, it proves that the solution provided is amenable to distributed implementation. This is of paramount significant as in an EV charging network the EVs are autonomous entities and operate based on the common information known to all the EVs (e.g., charging prices) and their local perception about network dynamics. Such distributed nature of the EV networks necessitates solving \eqref{mainProblem} in a distributed manner.

As it will be shown, an important feature of the proposed framework is that each EV adjusts its arrival rates solely based on the updated prices that are dynamically announced by the station operators, and such distributed adjustment of arrival rates by the EVs leads to an optimal global welfare for the entire network.

By recalling that utility function $U^n$ is  increasing in the arrival rates $\{\lambda_j\}$, decreasing in the LoLP $\beta_j$, concave in $\lambda_j$, and continuously differentiable in all of its arguments, the aggregate social welfare $R=\sum_{n=1}^N U^n$ is is maximized when $\forall n\in\{1,\dots,N\}$ and $\forall j\in\{1,\dots, J\}$
\begin{equation}\label{eq:derivative}
\frac{{\partial R}}{{\partial \lambda _j^{n} }} = \frac{{\partial {U^n}}}{{\partial \lambda _j^{n}}} + \sum\limits_{l = 1}^N \sum\limits_{s = 1}^J {\frac{{\partial {U^l}}}{{\partial {\beta_s}}}\cdot\frac{{\partial {\beta_s}}}{{\partial \lambda _j^{n}}} = 0} \ .
\end{equation}
Since the LoLP of each charge type depends only on the sum of the arrival rates of that type ${\lambda _j} = \sum\nolimits_{n = 1}^N \lambda  _j^{n}$, from \eqref{eq:derivative} we immediately have $\forall n\in\{1,\dots,N\}$ and $\forall j\in\{1,\dots, J\}$:
\begin{equation}\label{global}
\frac{{\partial R}}{{\partial \lambda _j^{n} }} = \frac{{\partial {U^n}}}{{\partial \lambda _j^{n}}} + \sum\limits_{l = 1}^N \sum\limits_{s = 1}^J {\frac{{\partial {U^l}}}{{\partial {\beta_s}}}\cdot\frac{{\partial {\beta_s}}}{{\partial \lambda _j}} = 0} \ .
\end{equation}
Solving \eqref{global} yields a {\em globally} optimal set of arrival rates for the customers.

\subsection{Local Problems}

We show that the globally optimal solutions yielded by \eqref{global} can be found in a distributed way, in which each charging statin forms and solves a local problem. In order the lay the foundation, we first formulate the local problems in this subsection, and relegate establishing the connection between the local and global solutions to Section III-C. In order to enforce such arrival rates the charging stations compute and announce prices for each customer type, which in turn forces the customers to operate at the desired (optimal) arrival rates. In order to solve this global problem and formalize such dynamic between posted prices and adjusted arrival  rates, by recalling that  $p_{j}$ is the prices charged to customers of type $j$, we first define the following {\em local} optimization problem for each customer $n\in\{1,\dots,N\}$ and service type $j\in\{1,\dots,J\}$:
\begin{align}\label{localOptim}
\tilde U^n({\boldsymbol \lambda}^{n};{\boldsymbol \beta}({\boldsymbol \lambda}))\;\triangleq\; \underbrace{U^n({\boldsymbol \lambda}^{n};{\boldsymbol \beta}({\boldsymbol \lambda}))}_{\text{gain}}\;-\;\underbrace{\sum\limits_{l = 1}^J {{p_{l}}\lambda _l^{n}(1 - {\beta_l}})}_{\text{cost}}\ .
\end{align}
By solving this local optimization problem, customer $n$ computes its locally optimal and relevant arrival rates ${\boldsymbol \lambda}^n=[\lambda_1^{n},\dots,\lambda_J^{n}]$. Specifically, the locally optimal arrival rates that maximize the {\em local} optimization problems $\{\tilde U^n\}_n$ for all $\forall n\in\{1,\dots,N\}$ and $j\in\{1,\dots,J\}$ satisfy:
\begin{equation}\label{localOptim2}
\frac{{\partial {U^n}}}{{\partial \lambda _j^{n}}} - {p_{j}}(1 - {\beta_j}) = 0\ .
\end{equation}
By solving \eqref{localOptim2} each customer $n$ can compute its own arrival rates based on the announced prices. The structure of the optimal local arrival rates depends on the selection of the utility function $U^n$, which is discussed in more details in Section~\ref{Example}.
\subsection{Connection Between Local and Global Solutions}
Before proceeding to the details of the relevance of the globally optimal solution to these local solutions it is important to note that individual customers do not have knowledge about the derivatives of their arrival rates with respect to LoLP $\beta_j$ which leads to above local optimization problem~\cite{courcoubetis1999pricing}. Also, the size of the system is fairly large when compared to demand of each individual customer, hence each EV does not have significant impact on the LoLP. Hence, by comparing equations (\ref{global}) and (\ref{localOptim2}) we can deduce that by the optimal prices for $j\in\{1,\dots,J\}$ satisfy:
\begin{equation}\label{prices}
p_{j}^{{*}} =  - (1 - {\beta^{j}})^{ - 1} \sum\limits_{l = 1}^N \sum\limits_{s = 1}^J {\frac{{\partial {U^l}}}{{\partial {\beta^{s}}}}\cdot\frac{{\partial {\beta^{s}}}}{{\partial \lambda _j}} = 0} \ .
\end{equation}
When a customer of type $j$ is presented with price $p_{j}^*$ the solution of her local optimization problem \eqref{localOptim} satisfies the conditions given in \eqref{eq:derivative}, which in turn guarantees that the solution of the global optimization problem in (\ref{global}) is equivalent to the local ones given in \eqref{localOptim} as established in~\cite{courcoubetis1999pricing}. Hence, the local optimal solutions become the equilibrium since no single customer can obtain a higher utility by deviating from its locally optimal solution. It is noteworthy that the prices $\{p_{1}^*,\dots,p_{J}^*\}$ can be considered as congestion prices since each customer pays other customers for one unit of marginal decrease in their utility due to the rise in the LoLP due to increase in the arrival rates.

\section{Resource Provisioning Framework}\label{sec:res}
\nomenclature{$Q^j_{C}$}{Number of simultaneous customers of type $j$ requesting $b_j$ units}%
\nomenclature{$S$}{Aggregated customer load given in \eqref{offeredLoad}}%
\nomenclature{$q_j$}{Average number of arrival rates of class $j$}%
In the resource provisioning problem characterize in \eqref{argmin} the optimal solution $C^*$ increases monotonically as the target outage upper bounds increase. Hence, one naive approach to solve this problem is a linear brute force search, which can be computationally prohibitive for large scale system. Instead, we aim to provide an establish an analytical connection between the optimal solution and the target outage requirements. This proposed approach, moreover, is versatile enough to be applied to other relevant optimization scenarios, e.g., integration of renewable generation or energy storage system and for varying customer demand.

In order to furnish the tools we model the request arrivals of the customers as a collection of $J$ independent queues (see~\cite{kleinrock}). In this model we start the analysis by first assuming that the grid resources are infinite and treat the case of finite resources in the next stage. Let $Q_\infty ^j$ denote the number of customers of type $j$ requesting $b_j$ units of power concurrently. Further, let $S$ be the sum of offered load by the system that is given by
\begin{equation}\label{offeredLoad}
S = \sum\limits_{j = 1}^J {{b_j}Q_\infty ^j}\ .
\end{equation}
Due to Poisson model for the arrival rates, for the mean and the variance of  $Q_\infty ^j$ we have $\mathbb{E}\left[ {Q_\infty ^j} \right]$ = ${\rm var}\left[ {Q_\infty ^j} \right]$ = ${q_j} = \frac{{{\lambda _j}}}{{{\mu _j}}}$. Therefore
\begin{equation}\label{stats}
\mathbb{E}\left[ S \right] = \sum\limits_{j = 1}^J {{b_j}{q_j}}, \hspace{2mm}{\rm{ and }}\hspace{3mm}{\rm var}\left[ S \right] = \sum\limits_{j = 1}^J {b_j^2} {q_j}.
\end{equation}
In the case of finite resources $C$, we can reformulate the carried load on the system by defining $Q_C^j$ as the number of simultaneous customers requesting $b_j$ units of grid resources. Then the blocking (loss-of-load) event occurs if upon an arrival of customer $j$ the load on the system  is greater than available resources $\sum\nolimits_{j = 1}^J {Q_C^j} $ is greater than $C-b_j$. For a given set of charging rates  $\boldsymbol{b}\triangleq [b_1,\dots, b_J]$ and average numbers of arrival rates $\boldsymbol{q}\triangleq [q_1,\dots, q_J]$,  let us denote the LoLP by of  customer type $j$ by $\beta_j:\mathds{R}^J \to [0,1]$, which captures the connection between the block probability on one hand hand, and  $\boldsymbol b$ and $\boldsymbol q$, on the other hand.

By noting that the random number of customers for each class $Q_\infty ^j$ are mutually independent Poisson random variables we have
\begin{align}
\beta _j(\boldsymbol{q},\boldsymbol{b}) &= {\mathds P}\left\{  {C - {b_j} < \sum_{j = 1}^J {{b_j}Q_C^j} } \right\} \label{blockings1}\\
&= {\mathds P}\left\{ {C - {b_j} < \sum_{j = 1}^J {{b_j} Q_\infty ^j} \le C} \;\Big|\; \sum_{j = 1}^J {{b_j}Q_\infty ^j \le C}\right\}\label{blockings2}\\
&=\frac{{{\mathds P}\left\{ {C - {b_j} < \sum_{j = 1}^J {{b_j}Q_\infty ^j}  \le C } \right\}}}{{{\mathds P}\left\{ {\sum_{j = 1}^J {{b_j}Q_\infty ^j \le C} } \right\}}}\label{blockings3} \ .
\end{align}

In the resource provisioning problem of interest, in order to meet multi class QoS targets $\{\delta_j\}$, $C^*$ should be at least as big as the mean offered load on the system, that is $S=\sum\nolimits_{j = 1}^J {{b_j}{q_j}} $. However, customers arrive in a stochastic fashion, hence it is required to add extra capacity to accommodate the fluctuations beyond the average offered load. For this purpose,  $C^*$ is set equal to the mean of the total load $\mathbb{E}\left[ S \right]$ adjusted by an extra term which is a multiple of its variance denoted by $x\cdot{\rm var}\left[ S \right]$. In this formulation more stringent QoS targets lead to larger values of  $x$. The objective of the resource provisioning problem is to provide a closed form solution to~\eqref{argmin}. To this end, we first scale the system with $\varsigma>0 $, according to which we have the following capacity for the network:
\begin{equation}\label{}
{\bar{C}(\varsigma,x)} = \varsigma \sum\limits_{j = 1}^J {{b_j}{q_j} + x\sqrt {\varsigma \sum\limits_{j = 1}^J {b_j^2{q_j}} } }\ .
\end{equation}
and we have the following limiting result \cite{ErlangMitra, hampshire}:
\begin{equation}\label{result1}
\lim_{\varsigma  \to \infty } \sqrt \varsigma  \beta _j(\boldsymbol{q},\boldsymbol{b}) = \frac{{{b_j}}}{{\sqrt {\sum\nolimits_{j = 1}^J {b_j^2{q_j}} } }}\cdot \frac{{\phi (x)}}{{\varphi (x)}}
\end{equation}
where $\phi (x) = \frac{1}{\sqrt{2\pi }}{e^{ -x^2/2}}$ and $\varphi (x) = \frac{1}{\sqrt{2\pi }}\int_{ - \infty}^x e^{-t^2/2}\; dt $. Now let us define function $\psi $ as the inverse of $\frac{\phi }{\varphi }$ for $\forall$$x$, that is:
\begin{equation}
\frac{{\phi (\psi (x))}}{{\varphi (\psi (x))}} = x \ .
\end{equation}
Note that $\psi(\cdot)$ is a strictly decreasing function with $\psi (y)+y>0$, $\forall y$ (\cite{ErlangMitra, hampshire } and references therein). Then, the asymptotic behavior of the QoS constraint $\beta^j(\boldsymbol{q},\boldsymbol{b})\leq\delta_j$ results:
\begin{equation}\label{xxx}
x \ge \psi \left( {\frac{{{\delta _j}}}{{{b_j}}}\sqrt {\sum\limits_{j = 1}^J {b_j^2{q_j}} } } \right) \ .
\end{equation}
The provisioned grid resources should satisfy QoS targets for all classes. Hence, the inequality in~\eqref{xxx} yields
\begin{equation}
x \ge \psi \left( {\mathop {\min }\limits_{1 \le j \le J} \frac{{{\delta _j}}}{{{b_j}}}\sqrt {\sum\limits_{j = 1}^J {b_j^2{q_j}} } } \right) \ .
\end{equation}
By operating at the lower possible values of $x$, the minimum required amount of resources solving the provisioning problem in~(\ref{argmin}) is
\begin{equation}\label{minCapacity}
{C^*} = \sum\limits_{j = 1}^J {{b_j}{\frac{{{\lambda _j}}}{{{\mu _j}}}} + \psi \left( {\mathop {\min }\limits_{1 \le j \le J} \frac{{{\delta _j}}}{{{b_j}}}\sqrt {\sum\limits_{j = 1}^J {b_j^2{\lambda _j}} } } \right)\sqrt {\sum\limits_{j = 1}^J {b_j^2{\lambda _j}} } } \ ,
\end{equation}
where $\psi(\cdot)$ can be easily computed numerically by solving~\cite{tian2007analysis}:
\begin{equation}
{x^{ - 1}}{e^{ - 0.5\psi {{(x)}^2}}} - \sqrt {2\pi }\;{\rm erf}\left( {\frac{1}{{\sqrt 2 }}\psi (x)} \right) - x\sqrt {0.5\pi }  = 0 \ ,
\end{equation}
and it can be plugged back into (\ref{minCapacity}) to compute the required capacity.

It is important to notice that providing enough resources to meet the QoS targets of the most dominant classes, that are the ones with minimum $\delta_j$/$b_j$, is sufficient for the remaining customer classes.

\section{Computing Loss-of-Load-Probabilities}\label{lossComp}

Both the pricing-based control problem (summarized in section \ref{pricing} with \eqref{mainProblem}) and the resource provisioning problem (summarized in \label{resourceProv} with with \eqref{argmin}) fall into multi-dimensional loss systems (or multi-rate Erlang loss systems), which are used to evaluate the guaranteed performance of networks with limited resources. In such systems an arriving customer requesting a certain amount of grid resources is either admitted to the system or is blocked and  we are interested in computing LoLP as a function of system parameters.

Computing the LoLP function requires the analysis of the $J$ independent time reversible Markov chains in which the state of the system is defined as the number of customers of each type, that is $\boldsymbol{Q}\triangleq [Q_C^1,\dots, Q_C^J ]$ and the state space is denoted by $\Omega\triangleq \{ {\boldsymbol{Q}: \sum_{j = 1}^J {{b_j}Q_\infty ^j \le C} } \}$. Let $\tilde{Q}_C^j$ denote the maximum number of customers of type $j$ that can be served simultaneously. Assuming the order $b_1\geq\dots \geq b_J\geq 0$, without loss of generality, provides that $0 \leq \tilde{Q}_C^1\leq  \dots  \leq \tilde{Q}_C^J$. Then the probability of being at state $\boldsymbol{Q}$ is~\cite{kleinrock}:
\begin{equation}
 {{\overline{\pi}}}(\mathbf{{Q}})=\prod\limits_{j = 1}^J {\frac{{q_j^{{Q_\infty^j}}}}{{{Q_\infty^j}!}}{e^{ - {q_j}}}} \ .
\end{equation}

Next, similar to (\ref{blockings2}) we condition on the finite capacity, and compute a generic state $\boldsymbol{Q}$ probability distribution as:
\begin{equation}
{ \pi}(\boldsymbol{Q})=\frac{{{\mathop{\overline{\pi}}\nolimits} (\mathbf{{Q}})}}{{\sum\nolimits_{ \tilde{Q} \in \Omega } { { \overline{\pi}}(\mathbf{\tilde Q})}}}  \cdot
\end{equation}
Next,  let us define the blocking states (LoLP) for customer type $j$ as
\begin{align*}
\Psi_j = \{{\boldsymbol{Q}: { {C - {b_j} \;<\; \sum\limits_{k= 1}^J {{b_k}Q_C ^k} \; \leq \; C}}}  \}\ .
\end{align*}
Hence, (\ref{blockings3}) can be re-written as:
\begin{equation}\label{block}
{\beta^j}(\boldsymbol{q},\boldsymbol{b}) = \sum_{s \in {\Psi _j}} {\pi (s)}  = 1 - \sum_{s \notin \Psi _j} {\pi (s)}\ ,
\end{equation}

where the second term represents the probability that the charging station has a total capacity of  $C-b_j$ instead of $C$ and $\pi(s)$ is the steady state probability mass function.
Furthermore, let us define function $H(C,J)$ as
\begin{equation}\label{eq:H}
H(C,J) \triangleq  \sum_{\left\{ {\boldsymbol{Q}:\;{\boldsymbol b}\boldsymbol{Q}\; \leq\; C} \right\}}{\prod\nolimits_{j = 1}^J {\frac{{q_j^{{Q^j}}}}{{{Q^j}!}}} } \centering\ ,
\end{equation}
based on which one can compute the LoLP for class $j$ explicitly as:
\begin{equation}\label{betaResult}
{\beta_j}(\boldsymbol{q},\boldsymbol{b})=1 - \frac{{H(C - {b_j},J)}}{{H(C,J)}} \cdot
\end{equation}
The set $\{ {\boldsymbol{Q}:\;{\boldsymbol b}\boldsymbol{Q}\; \leq\; C}\} $ in \eqref{eq:H} contains all the states corresponding to which {\em no} outage occurs.  While (\ref{betaResult}) provides an explicit representation for ${\beta_j}$, the associated computation can be costly when the system capacity $C$ is large. Considering a real-world scenario in which the number of classes $J$ is typically varies between 3 and 5, and $C$ is in the order of Mega-watts, computing LoLP could be easily carried out via the Kaufman-Roberts algorithm (Algorithm~\ref{algo}) which involves a simple recursion which computes the occupied resources~\cite{kaufman}. Let $c$ denote the amount of resources being in use and
\begin{align*}
\alpha(c)\;\triangleq\; {\mathds P}\left\{\mbox{c units of power in use} \right\}\ .
\end{align*}
We subsequently have
\begin{equation}\label{kauf}
\alpha(c)=\sum_{\left\{ {\boldsymbol{Q}:\;{\boldsymbol b}\boldsymbol{Q} \;\leq\; c} \right\}}{\frac{{q_j^{{Q^j}}}}{{{Q^j}\!}} \cdot \frac{1}{{H(C,J)}}}\ ,
\end{equation}
and the LoLP corresponding to customer type $j$ can be calculated by using
\begin{equation}
\beta_j(\boldsymbol{q},\boldsymbol{b}) = \sum\limits_{i = 0}^{{b_j} - 1} {\alpha (C - i)} \ .
\end{equation}
Note that the above derivations are based on the assumption that grid resources are discretized (e.g., $1$ kW is considered as $1000$ discrete serving units) and the interpretation of the Algorithm~\ref{algo} is that LoLP of customer type $j$ equals to some of occupied states that are in the range of $C-(b_j-1) ,\dots, C-1, C$.
Another important measure of interest for solving the optimal pricing problem given in \eqref{prices} is the set of derivatives of LoLP with respect to offered load by each customer class. In the proposed multi-class customer model, there is no explicit formulae for performance measures (LoLP) in terms of input parameters ($C$, $\lambda_{j}$,$\mu_{j}$, etc.). In this paper we follow the methods based on convolution algorithms~\cite{iversen2007derivatives} to compute the derivatives of the LoLP. To that end the derivative of the LoLP associated with customer type  $j$ with respect to the traffic intensity of another class $j_{1} \ne j_{2}$ can be computed by using the function $\alpha(\cdot)$ as follows.
\begin{eqnarray}\label{derivatives}
\frac{{\partial {\beta _{j_{1}}}}}{{\partial {q_{j_{2}}}}} &= \alpha (C - {b_{{j_1}}}) + \alpha (C - {b_{{j_{2}}}} - 1) +  \cdots  \\
&+ \alpha (C - {b_{{j_{2}}}} - {b_{j_{1}}} - 1) - (1 - {\beta ^{{j_{2}}}}){\beta ^{j_{1}}}\nonumber
\end{eqnarray}
A useful property of the LoLP is the elasticity property, which entails that the sensitivity of customer LoLP~$j$ to customer class $k$ is the same as the sensitivity of customer LoLP~$k $ to type $j$\cite{mazumdar}. This property is given as $\frac{{\partial {\beta_{{j_{1}}}}}}{{\partial {q_{{j_{2}}}}}} = \frac{{\partial {\beta_{{j_{2}}}}}}{{\partial {q_{{j_{1}}}}}}$.

\begin{algorithm}[t]\label{algo}
\begin{small}
\caption{Kaufman-Roberts Algorithm~\cite{kaufman}}\label{algo}
\begin{algorithmic}
\STATE Set $\kappa(0)=0$ and $\kappa(i)=0$ for $i \in {\rm I\!R}^{-}$
\FOR {$i$=$1$ to $C$}
\STATE $\kappa(i) = \frac{1}{i}\sum\nolimits_{j = 1}^J {{b_j}} {q_j}(j - {b_j})$
\ENDFOR
\STATE Compute $H = \sum\nolimits_{i = 1}^C {\kappa(i)} $
\FOR {$i$=$0$ to $J$}
\STATE $\alpha(i)=\frac{{\kappa (i)}}{H}$
\ENDFOR
\FOR {$j$=$1$ to $J$}
\STATE $\beta _j(\boldsymbol{q},\boldsymbol{b})=\sum\nolimits_{i = C - {b_j} + 1}^C {\alpha (i)} $
\ENDFOR
\end{algorithmic}
\end{small}
\end{algorithm}

Notice that Kaufman-Roberts algorithm is used in both of the proposed frameworks and the complexity of this algorithm is  \BigO{CJ} which provides a considerable improvement over solving the LoLP through \eqref{betaResult} which has a complexity of \BigO{C^J} ~\cite{kaufman,nilsson}. For the pricing-based control problem, we compute the derivatives of multi-rate blocking probabilities using the convolution-based algorithm \cite{iversen2007derivatives} and similar to Kaufman-Roberts algorithm, the order of complexity is \BigO{CJ}.

\section{Simulation Results}\label{Example}
 \begin{figure*}[t]
        \centering
                \begin{subfigure}[b]{0.32\textwidth}
                \centering
                                \includegraphics[width=\columnwidth]{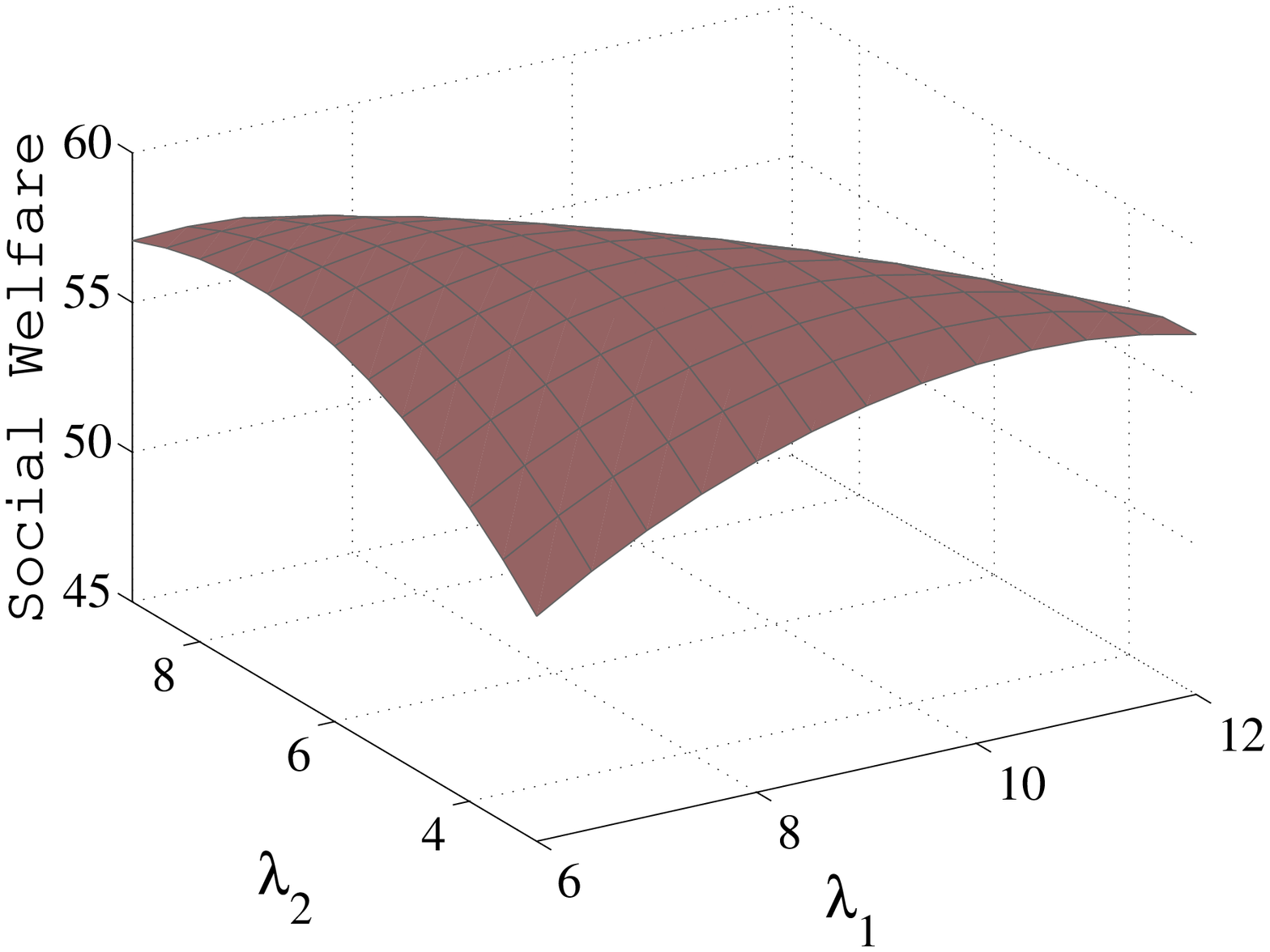}
 \caption{Social welfare computation.}\label{socialWelf}
                     \end{subfigure}
 \;
        \begin{subfigure}[b]{0.32\textwidth}
                \centering
                 \includegraphics[width=\columnwidth]{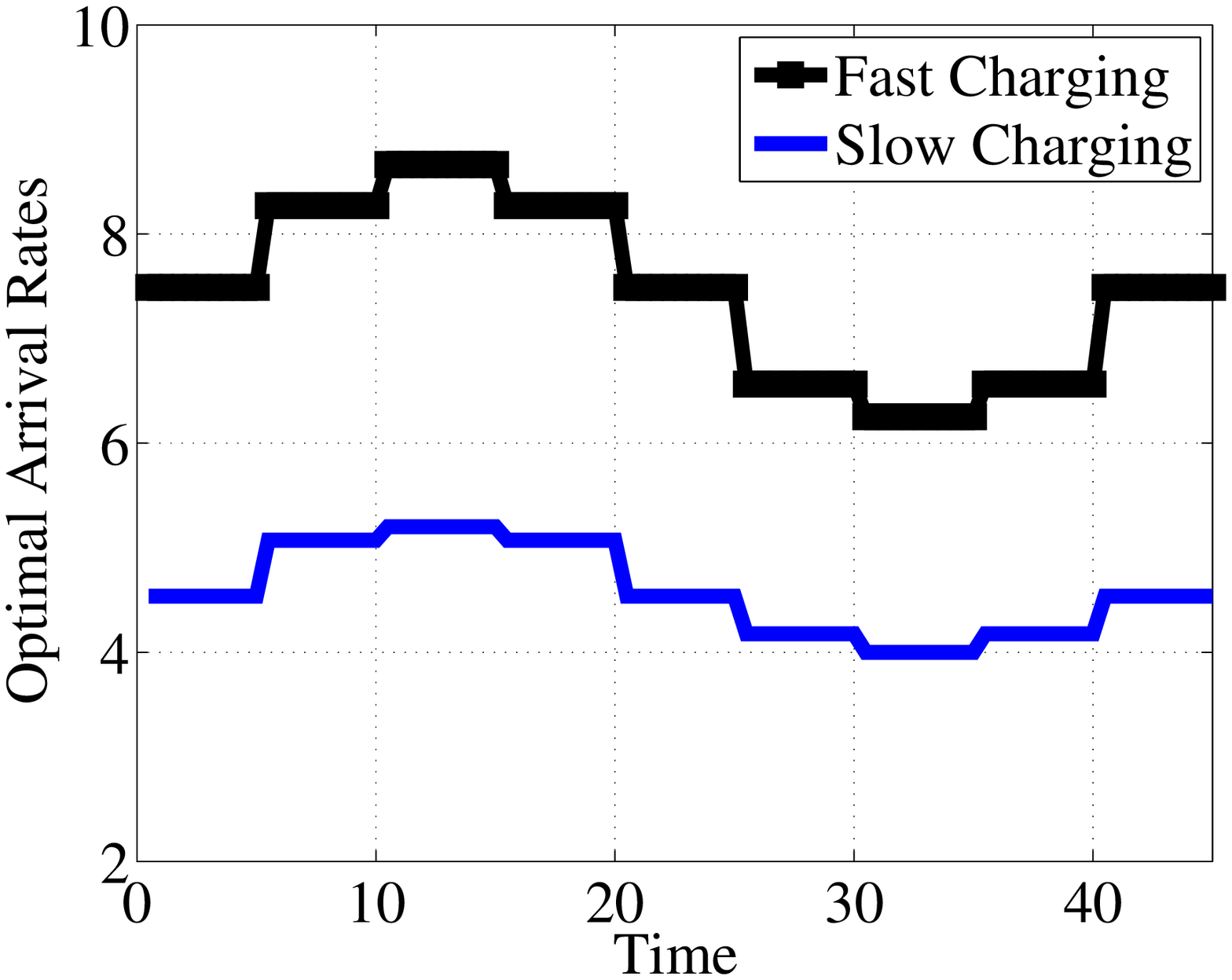}
 \subcaption{Optimal arrival rates ($\lambda_{1,2}$) with two classes.} \label{optRates}

        \end{subfigure}%
        \;
        \begin{subfigure}[b]{0.32\textwidth}
                \centering
                \includegraphics[width=\columnwidth]{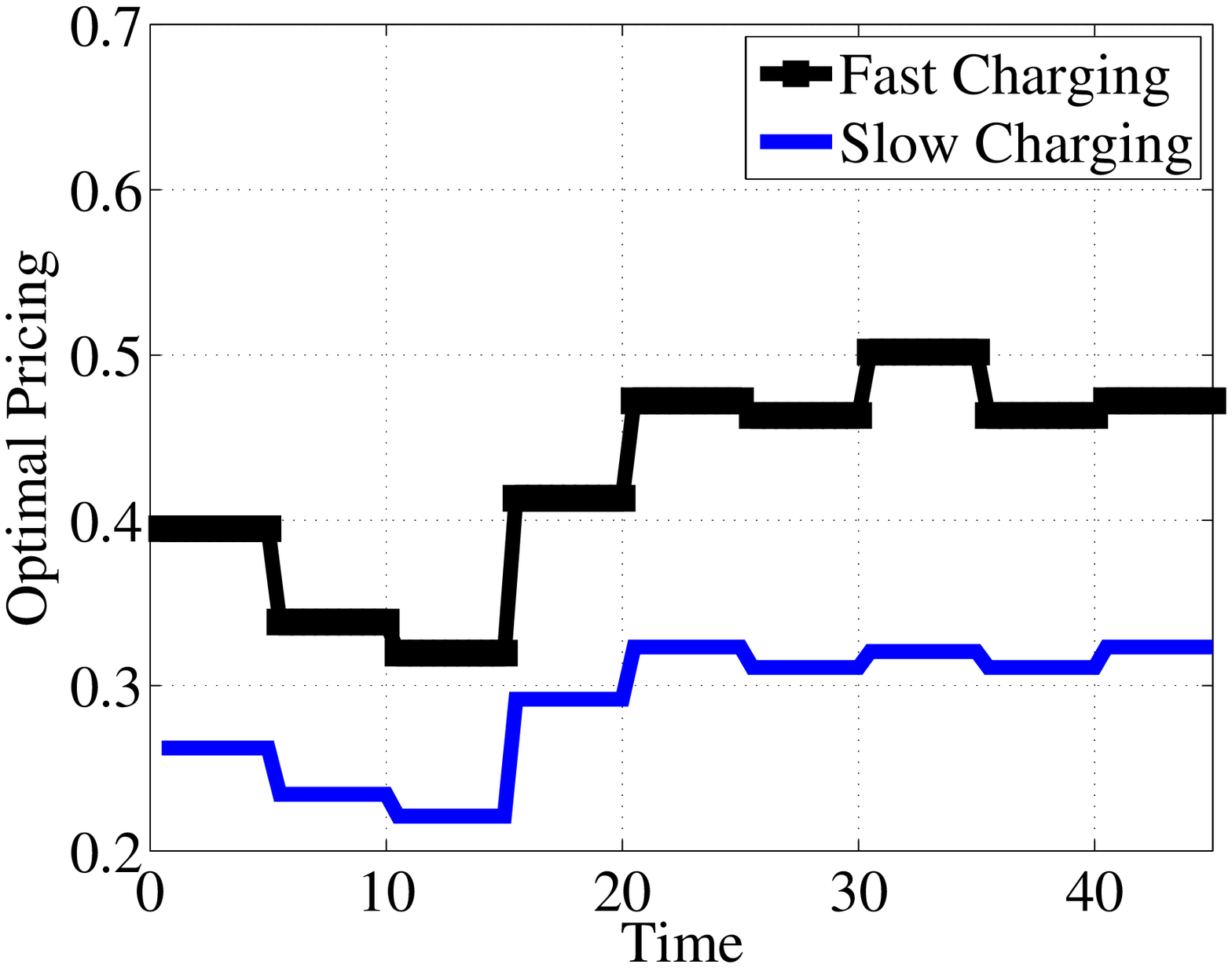}
 \subcaption{ Optimal arrival rates ($\lambda_{1,2}$) with two classes.}\label{pricing}
                     \end{subfigure}

        \caption{Numerical Evaluation-I: system capacity $C_{t\in k}$=$450+50\sin(2\pi \bar{k} / 80)$, $\bar{k}\in{(0,10,\dots,80)}$ and $C$ is assumed constant for a duration of $T$=$10$}\label{results1}
\end{figure*}

\subsection{Case Study-I: Pricing-Based Control}

In this subsection we present a case study to evaluate the social welfare maximization problem of the pricing-based control framework presented in section~\ref{sec:price}. The parameter setting for our case study is as follows. For the ease of representation we assume that there are two customer types, namely, fast charging customers (type-I) and slow charging customers (type-II). In order to mimic the typical charging rate for a fast DC charging (that is $50$kW), we set $b_1$=$50$ units. Also, charging duration takes around $20$min, so mean service rate is set as $\mu_1$=$3$. In a similar manner we tune the parameters for the slow charging customers as  $b_2$=$7$ units and $\mu_2$=$0.42$. For the utility functions we adopt the widely-used logarithmic utility function~\cite{fan2011, cong2}. The utility function of a single customer increases with the arrival rate and decreases with the customer LoLP according to
\begin{equation}
U=\left\{ {\begin{array}{ll}
{\sum\nolimits_{j = 1}^J {\omega_j\log (1 + {\lambda _j}) - {\theta _j}\log (1 + \beta _C^j)} }&{{\lambda _j} > 0}\\
&\\
0&{{\lambda _j} \le 0}
\end{array}} \right.\ ,
\end{equation}
where $\omega_j$ and $\theta_j$ are the weights of each class. Note that the weight is higher for customer types with higher demand $b_j$, mainly because fast charging customers demand more resources and hence gains more utility. A simple example is presented to clarify the matters. Assume that system capacity $C=500$, and the weights are chose as $\omega_1$=$20$, $\omega_2$=$10$ and $\theta_1$=$60$, $\theta_2$=$20$. Then the social welfare could be maximized by setting the arrival rates to $\lambda_1$=$8.6638$ and $\lambda_2$=$5.2001$. For this arrival rates prices $p_1$ and $p_2$ are computed to be $0.3197$ and $0.2211$ and the resulting LoLPs are $\beta^1$=$0.0097$ and $\beta^3$=$0.0009$, and the maximum utility is computed to be $59.1238$ units. Notice that $\omega_j$ is chosen greater than $\theta_j$ so that station operator is motivated to provide a good level of QoS. We further explore the relationship between the optimal arrival rates and the social welfare. As shown in Fig. \ref{socialWelf}, if the arrival rates deviate from its optimal value, the social welfare reduces.

 \begin{figure*}[t]
        \centering
                \begin{subfigure}[b]{0.32\textwidth}
                \centering
                   \includegraphics[width=\columnwidth]{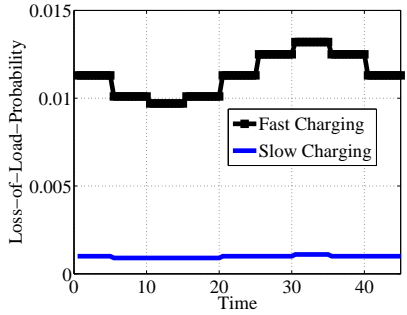}
 \caption{LOLP performance.}\label{blocking}

                      \end{subfigure}
\;
        \begin{subfigure}[b]{0.32\textwidth}
                \centering
\includegraphics[width=\columnwidth]{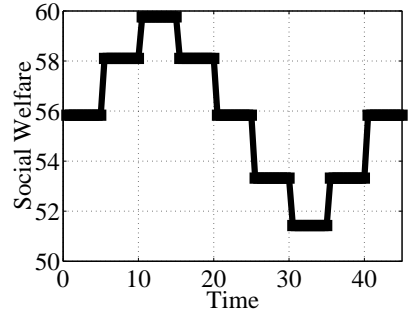}
 \caption{Social welfare computation.}\label{TotalWelfare}

        \end{subfigure}%
        \;
        \begin{subfigure}[b]{0.32\textwidth}
                \centering
                \includegraphics[width=\columnwidth]{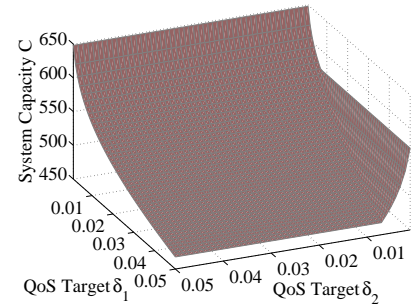}
         \caption{Capacity planning for varying QoS targets with two classes.}\label{capPlanning}
        \end{subfigure}
        \caption{Numerical Evaluation-II}\label{results1}
\end{figure*}

We proceed to provide more numerical evaluations. In the first setting we present the relationship between system capacity and optimal arrival rates for time-dependent case, where station capacity varies over time and follows $C(k)$=$450+50\sin(2\pi k / 80)$ (due to grid conditions) and also $C(k)$ is assumed to be constant for every $T=10$ duration. Results depicted in Fig. \ref{optRates} shows that, it is more beneficial to accept fast charging customers as they improve the social welfare function more than slow charging ones. For the given set of arrival rates, the corresponding prices in \eqref{prices} and loss-of-load-probabilities are given in Figs. \ref{pricing} and \ref{blocking} respectively. Obviously, since fast charging uses more resources, the corresponding prices are higher and due to higher arrival rate the corresponding LoLP is higher than the slow charging customers. Finally, we present the corresponding social welfare in Fig. \ref{TotalWelfare}.

\subsection{Case Study-II: Capacity Planning}

We proceed to compute minimum amount of grid resources to provide QoS guarantees to two customer classes with the same set of parameters ($\mu_1$=$3$, $\mu_2$=$0.42$, $b_1$=$50$, $b_2$=$7$) for a wide range of QoS targets ($0.001 \le {\delta _1} \le 0.05$, $0.001 \le {\delta _2} \le 0.05$) with fixed arrival rate $\lambda_{1,2}$=$5$. The results depicted in Fig.~\ref{capPlanning} presents that, since type-$1$ is the dominant class in (most) regions where $\delta_1/b_1<\delta_2/b_2$, providing resources for  class-$1$ already satisfies the QoS targets for class-$2$. Next we compute the required capacity for different range of arrival rates and fixed QoS targets $\delta_1$=$\delta_2$=$0.03$. Results depicted in Fig.~\ref{capPlanningX} can be used as a guideline to choose the required capacity for given arrival rates.

We proceed to investigate the percentage of reduction in station capacity ($C$) for a given LoLP targets. The motivation is that instead of providing zero percent {LoLP, station operator can sacrifice to reject small amount of customers and reduce the stress on the grid. We use the same parameter setting as above and computed the required capacity respect to almost zero LoLP ($\delta_{1,2}$=$10^{-6})$. As presented in Fig.~\ref{capPlanning2} even one percent QoS targets leads to significant savings in station capacity. Our final evaluation is on evaluating the system performance for non-homogenous arrival rates. Let's assume that arrival rate for class-II is constant and $\lambda_2(k)$=$10$ and  $\lambda_1$=$10+2\sin(2\pi k / 80)$. The QoS targets are set as $\delta_1$=$0.04$ and $\delta_2$=$0.01$. Then station operator can provision the system according to peak hour demand (computed as $C$=$683$) and hence meet the QoS requirements at all times. Results are presented in Fig.~\ref{lambdas}.


\begin{figure*}[t]
\centering
 \begin{subfigure}[c]{0.32\textwidth}
                \centering
\includegraphics[width=\columnwidth]{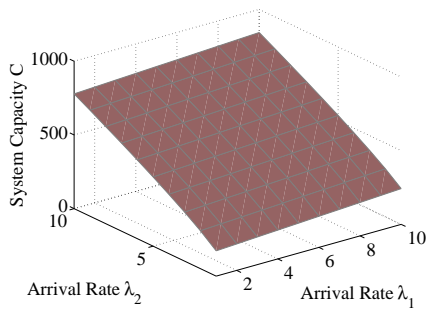}
 \caption{Resource provisioning (fixed $\delta_{1}=\delta_{2}=0.03$)}\label{capPlanningX}
   \end{subfigure}
        \begin{subfigure}[c]{0.32\textwidth}
                \centering
                \includegraphics[width=\columnwidth]{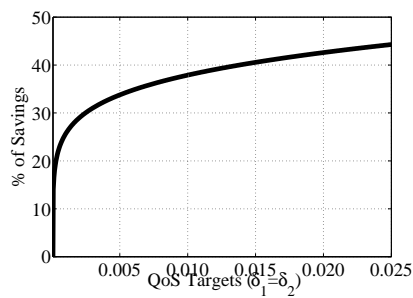}
 \caption{\% of Savings in Station Capacity}\label{capPlanning2}
        \end{subfigure}%
        \begin{subfigure}[c]{0.32\textwidth}
                \centering
                \includegraphics[width=\columnwidth]{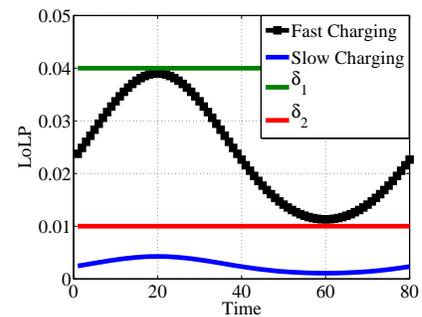}
 \caption{LoLP Performance, $C_{min}=683$.}\label{lambdas}
        \end{subfigure}
        \caption{Numerical Evaluation III}\label{results1}
\end{figure*}

\section{Concluding Remarks}
In this paper we provided two important design problems for electric vehicle charging stations with multiple classes of customers. We considered two cases. In the first one, given a capacity of a station and infinite amount of customer request, we provide a framework to compute the optimal arrival rates such that the total social welfare is maximized. This is a very typical case for charging stations located in big cities. On the other hand, in the second case we considered charging stations located in small cities. This time our primary concern was to calculate the minimum amount of grid resources such that each customer class is ensured with a certain level of QoS targets (LoLP).

This initial work can be expanded in different directions. For example each station can employ an energy storage system that can aid to further reduce the strain on the power grid. The ESS can be charged during light traffic and the stored energy can be used to meet EV demand during peak hours. Another future research interest would be to consider a different resource policy to optimize the resource usage (e.g., portioning for each class etc.). In this paper we assume one micro-grid level charging. Another research direction can include a customer routing among the micro grids or stations.

In this work, we have concentrated our focus on a network of stations fed by a single substation. To that end, our final research direction is on considering a general network case and address the congestion issues in a grid composed of power lines and other elements with different capacity ratings.

\section*{Acknowledgment}
This publication was made possible by NPRP grant \# 6-149-2-058 from the Qatar National Research Fund (a member of Qatar Foundation). The statements made herein are solely the responsibility of the authors.

The authors would like to thank Dr. Sercan Teleke for the fruitful discussions.

\ifCLASSOPTIONcaptionsoff
  \newpage
\fi

\bibliographystyle{IEEEtran}
\bibliography{energycon}

\end{document}